\documentclass[12pt]{article} 

\newtheorem{theoA}{Theorem}
\newtheorem{theoB}{Theorem}

\newtheorem{predl}{Statement}
\newtheorem{prop}{Proposition}
\newtheorem{concl}{Corollary}
\newtheorem{lemma}{Lemma}
\newtheorem{note}{Remark}

\newcommand{\thRef}[1]{\ref{#1}}
\newcommand{\Ref}[1]{(\ref{#1})}

\newenvironment{proof}{\textbf{Proof. }}
{$\bigtriangleup$}
\newenvironment{eq}{\begin{equation}}{\end{equation}}

\newcommand{\si}{\sigma}
\newcommand{\Si}{\sigma_2}
\newcommand{\al}{\alpha}
\newcommand{\la}{\lambda}

\newcommand{\LA}{\langle}
\newcommand{\RA}{\rangle}
\newcommand{\ov}[1]{\overline{#1}}
\newcommand{\id}[1]{{{\rm id}\{{#1}\}}}

\newcommand{\tr}{{\rm tr}}
\newcommand{\sys}{{\mathcal{S}} }
\newcommand{\projlim}{\mathop{\rm{projlim }}}
\newcommand{\alg}{\mathop{\rm alg}}
\newcommand{\mdeg}{\mathop{\rm mdeg}}
\newcommand{\diag}{\mathop{\rm diag}}
\newcommand{\Char}{\mathop{\rm char}}
\newcommand{\sign}{\mathop{\rm{sgn }}}
\newcommand{\mod}{\mathop{\rm{mod }}}

\begin{document}
 \title{The Algebra of Invariants of $3\times3$ Matrices over a
 Field of Arbitrary Characteristic}
 \author{A. A. Lopatin\thanks{Supported by RFFI (grant 01.01.00674).}\\ Chair of Algebra,
 \\ Department of Mathematics, \\
 Omsk State University,\\ 55A Prospect Mira, Omsk 644077, Russia,
 \\e-mail: lopatin@math.omsu.omskreg.ru }

 \date{}

\maketitle

\begin{abstract} The least upper bound on degrees of elements of a minimal system of
generators of the algebra of invariants of $3\times3$ matrices is
found, and the nilpotency degree of a relatively free finitely
generated algebra with the identity $x^3=0$ is established.
\end{abstract}

\section{Introduction}
Let $K$ be an infinite field of arbitrary characteristic. Let a
reductive algebraic group $G$ act regularly on $m$-dimensional
affine variety $V=K^m$. This action defines natural action of $G$
on the coordinate algebra $K[V]$: $(g\cdot f)(v)=f(g^{-1}\cdot
v)$, where $ f\in K[V]$, $g\in G$, $v\in V$. Denote by $R=K[V]^G$
the algebra of invariants of $K[V]$ with respect to the action of
$G$. By the Hilbert--Nagata Theorem, it is a finitely generated
graded subalgebra. But Hilbert's proof for the case $\Char(K)=0$,
as well as Nagata's proof for the case $\Char(K)>0$, are not
constructive. The goal of the constructive theory of invariants is
to find a minimal (i.e. irreducible) homogeneous system of
generators (MSG) of $K[V]^G$ explicitly. It is an important
problem, which arose as early as the theory of invariants itself.
If one knows generators for each homogeneous component of the
algebra of invariants, then, theoretically, the problem of finding
MSG is equivalent to finding a constant $D$ such that $K[V]^G$ is
generated by invariants of degree not greater than $D$~\cite{Pop}.
Popov gave a bound $D$ for a connected semisimple group in
characteristic zero case~\cite{Pop}. But Popov's bound is rather
rough, so the problem of finding finer bounds is open.

Let $N_0=\{0,1,2,\ldots\}$. If $A$ is a $N_0$-graded algebra,
denote by $A^{+}$ the subalgebra generated by elements of $A$ of
positive degree. It is easy to see that the set $\{r_i\} \subseteq
R$ is a MSG iff $\{\ov{r_i}\}$ is a basis of
$\ov{R}={R}/{(R^{+})^2}$. Call an element $r\in R$ {\it
decomposable} if it belongs to the ideal $(R^{+})^2$. So the least
upper bound for the degrees of elements of MSG of the algebra of
invariants is equal to the highest degree of indecomposable
invariants.

Let $G=GL_n(K)$ act on the affine space
$M_{n,d}(K)=M_n(K)\oplus\cdots\oplus M_n(K)$ ($d$ times)  by
conjugations according to the following rule:
$B\cdot(A_1,\ldots,A_d)\rightarrow(BA_1B^{-1},\ldots,BA_dB^{-1})$,
where $M_n(K)$ is the space of all $n\times n$ matrices over $K$,
$A_i\in M_n(K)$, $B\in GL_n(K)$ $(i=\ov{1,d})$. This action
induces an action on the coordinate ring $K_{n,d}=K[x_{ij}(r)\mid
i,j=\ov{1,n} ;\; r=\ov{1,d}]$. Denote by $R_{n,d}=K_{n,d}^G$ the
algebra of invariants. Let $X_r=(x_{ij}(r))_{1\leq i,j\leq n}$ be
the generic matrices of order $n$, and let $\si_k(X)$ be the
coefficients of the characteristic polynomial of an $n\times n$
matrix $X$
$$\chi_n(X)=X^n-\si_1(X)X^{n-1}+\cdots+(-1)^n\si_n(X)E.$$
It is easy to see that $\si_k(X_{i_1}\cdots X_{i_s })\in R_{n,d}$.
Denote by $D(n,d,K)$ (by $D_{\si_k}(n,d,K)$, respectively) the
highest degree of indecomposable invariants (of indecomposable
invariants of the form $\si_k(X_{i_1} \cdots X_{i_s})$,
respectively).

In the case $\Char(K)=0$, $R_{n,d}$ has been well investigated.
Its generators and relations are described in~\cite{Sib},
\cite{Pro}, \cite{Raz}. In particular, it is known that $R_{n,d}$
is generated by its elements of degree $\leq n^2$~\cite{Raz}.
Great progress in the study of $R_{n,d}$ in the case of positive
$\Char(K)$ was made due to Donkin and Zubkov. Donkin showed that
$R_{n,d}$ is generated by all elements of the form
$\si_k(X_{i_1}\cdots X_{i_s })$~\cite{Don}, and Zubkov~\cite{Zub1}
extended Procesi-Razmyslov's Theorem on the relations to this
case. Before formulating this theorem, we must fix some notation.

Let $A$ be an associative algebra. Denote by
$A\mbox{-}\alg\{b_1,\ldots,b_d\}$ the associative algebra
generated over $A$ by $1,b_1,\ldots,b_d$ which commute with $A$.
If $b_1,\ldots,b_d$ are free over $A$, i.e. all identities of
$C=A\mbox{-}\alg\{b_1,\ldots,b_d\}$ are consequences of identities
of $A$, identities  $b_ia=ab_i$ ($i=\ov{1,d}$, $a\in A$) and
identities of associativity, we denote $C$ by $A\LA
x_1,\ldots,x_d\RA$. If $C=A\mbox{-}\alg\{b_1,\ldots,b_d\}$ is a
$N_0$-graded algebra with a unit over a $N_0$-graded subalgebra
$A$, and elements $b_1,\ldots,b_d$ are homogeneous of positive
degree, we denote by $C^{\#}$ the graded subalgebra $\sum_{1\leq
i\leq d}Cb_i$. The ideal generated by $a_1,\ldots,a_d$ is denoted
by ${\rm id }\{a_1,\ldots,a_d\}$. We call an identity $h=0$ of $A$
a {\it consequence} of identities $\{h_i=0\}$ if $h=\sum_{j=1}^s
a_jh_{k_j}a_j'$, where $a_j,a_j'\in A$. The homogeneous component
of degree $r$ of a graded algebra $A$ is denoted by $A(r)$.

The algebra of {\it concomitants} $C_{n,d}$ for $M_{n,d}$ is the
algebra of all polynomial $GL_{n}(K)$--equivariant mappings of the
space $M_{n,d}(K)$ to $M_{n}(K)$, where $GL_{n}(K)$ acts on
$M_{n,d}(K)$ and $M_{n}(K)$ by conjugation. It is easy to see that
$C_{n,d}$ is isomorphic to
$R_{n,d}\mbox{-}\alg\{X_1,\ldots,X_d\}$, i.e. the subalgebra of
$M_{n}(K_{n,d})$ generated by the generic matrices
$X_1,\ldots,X_d$ over $R_{n,d}$~\cite{Zub3}. For $k>n$, consider
the embedding $M_n\rightarrow M_k$ taking $A\in M_n$ to the matrix
whose left upper $n\times n$ cell coincides with $A$ and all the
other entries are equal to zero. This mapping induces
homomorphisms $R_{k,d}(r)\rightarrow R_{n,d}(r)$ and
$C_{k,d}(r)\rightarrow C_{n,d}(r)$ (consult~\cite{Zub3} for
details). Taking projective limits, we obtain the {\it free
algebra of invariants} of $d$ matrices
$R_d=\bigoplus_{r\geq0}{\projlim_n R_{n,d}(r)}$ and the {\it free
algebra of concomitants} $C_d=\bigoplus_{r\geq0}{\projlim_n
C_{n,d}(r)}$. Let $S$ be the free semigroup generated by letters
$\{a_1,a_2,\ldots\}$. Words $b=a_{i_1}\cdots a_{i_l}$ and
$c=a_{j_1}\cdots a_{j_l}$ are called equivalent, if there exists a
cyclic permutation $\pi\in S_l$ such that $i_s=j_{\pi(s)}$,
$s=\ov{1,l}$. The {\it cycle} (in letters $a_1,a_2,\ldots$) is the
equivalence class of some word. The cycle is {\it primitive}, if
it is not equal to a power of a shorter cycle. It is known that
$R_d$ is isomorphic to the algebra of polynomials in 'symbolic'
free generators $\si_k(a)$, where $a$ is a primitive cycle in
letters $x_1,\ldots,x_d$, $k>0$~\cite{Don2}. The algebra $C_d$ is
isomorphic to the free associative algebra in 'formal matrix
variables' $x_1,\ldots,x_d$ over $R_d$, i.e. $C_d\simeq R_d\LA
x_1,\ldots,x_d\RA$~\cite{Zub3}. If $n\geq r$, then $R_{n,d}(r)$ is
naturally isomorphic to $R_{d}(r)$, and $C_{n,d}(r)$ is isomorphic
to $C_{d}(r)$ (see~\cite{Don2}). Denote by $\deg(c)$ the degree of
a word $c$, i.e. the number of letters appearing in $c$, and by
$\mdeg(c)$ the multidegree of $c$, i.e.
$\mdeg(c)=(\lambda_1,\lambda_2,\ldots)$, where $\lambda_j$ is the
number of times $a_j$ appears in $c$. These notations are also
used for cycles.

Let $B$ be an arbitrary commutative algebra, $q_1,\ldots,q_s\in
B$, and let $A_1,\ldots,A_s$ be $n\times n$ matrices over $B$. For
$k=\ov{1,n}$ Amitsur's formula states~\cite{A}:
\begin{eq}\label{amiss}
\si_k(\sum_{l=1}^s q_l A_l)=\sum (-1)^{k-(j_1+\cdots+j_t)}q^{j_1
\mdeg(c_1)+\cdots+j_t\mdeg(c_t)}
\si_{j_1}(c_1)\cdots\si_{j_t}(c_t),
\end{eq}
where $q^{(\lambda_1,\ldots,\lambda_s)}=q_1^{\lambda_1}\cdots
q_s^{\lambda_s}$, and the sum ranges over all pairwise different
primitive cycles $c_1,\ldots,c_t$ in letters $A_1,\ldots,A_s$ and
numbers $j_1,\ldots,j_t$ with $\sum_{i=1}^{t}j_i\deg(c_i)=k$.
By~\Ref{amiss} one can express $\si_k(G)\in R_{n,d}$ in terms of
elements of the form $\si_k(U)$, where $G\in C_{n,d}^{\#}$, $U$ is
a non-empty word in the generic matrices. Identification $R_n(r)$
with $R_{n,d}(r)$ for $n\geq r$ allows one to define
$\sigma_k(g)\in R_d$ for $g\in C_{d}^{\#}$, $k\geq1$, correctly.

The algebra $R_d$ can be regarded as an associative-commutative
$K$--algebra with a unit generated by 'symbolic' elements
$\si_k(g)$, $g\in K\LA x_1,\ldots,x_d\RA^{\#}$, $k>0$. The ideal
of relations of the algebra $R_d$ is generated by
(see~\cite{Zub2}):

$(A)$ $\forall k\geq1, \forall g, h\in K\LA
x_1,\ldots,x_d\RA^{\#}, \si_k(gh)=\si_k(hg)$.

$(B)$ Amitsur's formula.

$(C)$ $\forall \alpha\in K, \forall k\geq1, \forall g\in K\LA x_1,
\ldots,x_d\RA^{\#}, \si_k(\alpha g)=\alpha^k\si_k(g).$

$(D)$ $\forall t,k\geq1,  \forall g\in K
    \LA x_1,\ldots,x_d\RA^{\#},
    \si_k(g^t)=\sum\limits_{i_1,\ldots,i_{kt}}\beta^{(k,t)}_{i_1,\ldots,i_{kt}}
    \si_1(g)^{i_1}\cdots\si_{kt}(g)^{i_{kt}}$, where coefficients
$\beta^{(k,t)}_{i_1,\ldots,i_{kt}}
    \in Z$ are determined uniquely.

The kernels of natural projections $R_d\rightarrow R_{n,d}$,
$C_d\rightarrow C_{n,d}$ we denote by $I_{n,d}$, $J_{n,d}$,
respectively. Procesi--Razmyslov's Theorem asserts that the ideal
$I_{n,d}$ is generated by 'symbolic' elements $\si_k(f)$, $k>n$,
and the ideal $J_{n,d}$ --- by elements
$\chi_k(f)=f^k-\si_1(f)f^{k-1}+\cdots+(-1)^k\si_k(f)$, $k\geq n$,
where $f\in C_{d}^{\#}$. In other words, the ideal of relations of
$R_{n,d}$ is generated by $(A)$--$(E)$, where

$(E)$ $\forall k>n, \forall f\in C_{d}^{\#}, \si_k(f)=0$.

In~\cite{Dom_K_Z} it is proved that, in contrast to the case 
of $\Char(K)=0$, if $0<\Char(K)\leq n$ then the degree bound for
the generators of $R_{n,d}$ tends to infinity when $d$ tends to
infinity. In~\cite{Dom_K_Z} an explicit MSG for $R_{2,d}$ is
given. In~\cite{Dom} some upper and lower bounds on $D(n,d,K)$ are
pointed out.

In this paper we consider the case $n=3$. We find the least upper
bound on degrees of elements of MSG of $R_{3,d}$ (except for the
case of $\Char(K)=3$, $d=6k+1$, $k>0$, where we estimate the least
upper bound with error not greater than $1$).

\begin{theoB} The least upper bound $D=D(3,d,K)$ on degrees of
elements of a minimal system of generators of the algebra of
invariants $R_{3,d}$ for $d>1$ is equal to:
\begin{enumerate}
\item[] if $\Char(K)=0$ or $\Char(K)>3$, then $D=6$,

\item[] if $\Char(K)=2$, then $D=\left\{
\begin{array}{ccc}
d+2&,&d\geq4\\
6&,&d=2\;{\rm or}\;3\\
\end{array} \right.$

\item[] if $\Char(K)=3$, then $D=\left\{
\begin{array}{ccl}
3d&,&d\;{\rm even}\\
3d-1&,&d\equiv 3\;{\rm or}\;5\;(\mod 6) \\
3d-1 \;{\rm or}\; 3d&,& d\equiv 1\;(\mod 6). \\
\end{array} \right.$
\end{enumerate}
If $d=1$, then $D=3$.

\end{theoB}


The {\it nilpotency degree} of a graded algebra
$A=\bigoplus_{j\geq0}A(j)$, where $A(0)=K$, is the least $C$ for
which $a_1\cdots a_C=0$ for all $a_i\in\ A^{+}$ ($i=\ov{1,C}$).
The idea of the proof of Theorem 2 consists in reduction of the
problem of decomposability of certain invariants to the problem of
equality to zero of certain elements of the algebra $N_{n,d} =
{K\LA x_1,\ldots,x_d \RA}/{\id{\;f^n\;|\; f\in K\LA
x_1,\ldots,x_d\RA^{\#}}}$ and, in particular, to the question of
finding $C(n,d,K)$ --- the nilpotency degree of $N_{n,d}$ (see
Lemmas ~\ref{cr_1},~\ref{cr_2},~\ref{cr_3}). In the case of
characteristic zero $n(n+1)/2\leq C(n,d,K)\leq n^2$
(see~\cite{Kuzmin},~\cite{Raz_book}), and Kuzmin conjectured that
$C(n,d,K)=n(n+1)/2$. This conjecture has been proved to be true
for $n\leq4$~\cite{Vau}. For prime characteristic, there exists an
upper bound on $C(n,d,K)$~\cite{Klein}: $C(n,d,K)<(1/6)n^6d^n$.

In Section 2 of this paper we prove the following theorem.
\begin{theoA}The nilpotency degree $C=C(3,d,K)$ of $N_{3,d}$
($d>1$) equals:
\begin{enumerate}
\item[]if $\Char(K)=0$ or $\Char(K)>3$, then $C=6$,
\item[]if $\Char(K)=2$, then
$C=\left\{
\begin{array}{lcl}
d+3&,& d\geq3\\
6&,& d=2,\\
\end{array}
\right.$ \item[]if $\Char(K)=3$, then $C=\left\{
\begin{array}{lcl}
3d+1&,& d{\rm\ is\ even}\\
3d\;{\rm or}\;3d+1 &,& d{\rm\ is\ odd}.\\
\end{array} \right.
$
\end{enumerate}
\end{theoA}

These theorems show one more difference between the cases of
characteristic zero and prime: for $\Char(K)=0$,
$D(3,d,K)=C(3,d,K)$ ($d>1$), while for $\Char(K)=2,3$
$D(3,d,K)<C(3,d,K)$ ($d>1$, and $d\neq2,3$ if $\Char(K)=2$).


\section{Associative algebra with the identity
$a^3=0$}\label{chapter_nilp}

\subsection{General remarks}\label{intro}
In this section we compute the nilpotency degree of a relatively
free finitely generated algebra $N_{3,d}={K\langle x_1, \ldots,
x_d\rangle}/{I}$, where $I=\mbox{id}\{f^3\mid\ f\in K\LA
x_1,\ldots,x_d\RA^{\#}\}$, satisfying 
$$T_1(a)=a^3=0,\;a\in K\LA x_1,\ldots,x_d\RA^{\#}.$$

We call $x_i$ letters, and monomials in $x_i$ words. By $x,y,z$ we
denote any triple of pairwise distinct letters. Throughout this
section, all considered elements of $N_{3,d}$ are meant to be
non-empty words, and all words are meant to belong to $N_{3,d}$,
if we do not explicitly write otherwise. Small Greek letters
(possibly with index) denote elements of $K$. Denote by $p$
($p=0,2,3,\ldots$) the characteristic of the field $K$.

Since the ideal $I$ is homogeneous, $N_{3,d}$ possesses natural
$N_0$- and $N_0^d$-gradings by degrees and multidegrees,
respectively, for which we use the same notations as in
Introduction. The degree of a word $u$ in a letter $x$ we denote
by $\deg_x(u)$. The multidegree ${(\al,\ldots,\al)}$ ($d$ times)
will also be denoted by $\al^{(d)}$.

Partial and complete linearization of $a^3=0$ gives the identities %
$$\begin{array}{c}
T_2(a,b)=a^2b+aba+ba^2=0.\\
T_3(a,b,c)=abc+acb+bac+bca+cab+cba=0.
\end{array}$$

Denote by ${\sys}$ the system
$$\left\{
\begin{array}{ccl}
g_1T_1(f)g_2&=&0\\
g_1T_2(f_1,f_2)g_2&=&0\\
g_1T_3(f_1,f_2,f_3)g_2&=&0,
\end{array}\right.\eqno(\sys)
$$
where $f$, $f_i$ are non-empty words $(i=1,2,3)$, words $g_1$,
$g_2$ can be empty,  and equalities are meant to hold modulo ideal
$I$. Let $\sys_{\Lambda}$ be the subsystem of $\sys$ which
consists of equations of multidegree $\Lambda$. For each word $u$
of multidegree $\Lambda$, introduce a variable $x_u$, and regard
system $\sys_{\Lambda}$ as a homogeneous system of linear
equations in $\{x_u\}$ over $K$. Clearly, if $\mdeg(u)=\Lambda$,
then $u=0$ in $N_{3,d}$ iff $x_u=0$ for each solution of
$\sys_{\Lambda}$. If $h=0$ is an equation from $\sys_{\Lambda}$,
by $h|_{\{a_u\}}$ we denote the result of substitution
$\{x_u=a_u\}$ in $h=0$, where $a_u\in K$.

We call a word {\it canonical} with respect to $x_i$, if it has
one of the following forms: $w_1$, $w_1x_iw_2$, $w_1x_i^2w_2$,
$w_1x_i^2ux_iw_2$, where subwords $w_1$, $u$, $w_2$ do not contain
$x_i$, and subwords $w_1$, $w_2$ can be empty. If a word is
canonical with respect to each $x_i$, we call it {\it canonical}.
\begin{predl}\label{canon}
An arbitrary word $w\in N_{3,d}$ is equal to a sum of canonical
words which belong to the same homogeneous component as $w$.
\end{predl}
\begin{proof}
$T_2(x_i,w)=0$ implies
\begin{eq}\label{a1}x_iwx_i=-x_i^2w-wx_i^2.\end{eq}
Since $T_2(x,xw)=0$, it follows that
\begin{eq}\label{a}xwx^2=-x^2wx.\end{eq}
Applying to each letter~\Ref{a1} and then~\Ref{a}, we obtain the
required. 
\end{proof}
\begin{note}
The presentation of a word from Statement~\ref{canon} does not
have to be unique.
\end{note}

\begin{concl}
If a word $w \in N_{3,d}$ contains more than $4$ occurrences of
some letter, then $w=0$. In particular, the length of a non-zero
word does not exceed $3d$.
\end{concl}

Hereafter, to specify the subword to which the identity is
applied, we sometimes put it in parentheses. Also, if we need to
split a word into a product of subwords, we insert dots in it.
(For example, see the deduction of~\Ref{5.1} from~\Ref{f3}.)
Moreover, we will apply Statement~\ref{canon} to all words without
reference.

Let us obtain some identities. Applying~\Ref{a1}, we get
\begin{eq}\label{1p1}
(xy)^2=(xyx)y=-x^2y^2-(yx^2y)=y^2x^2.
\end{eq}
Further, we apply~\Ref{1p1} to all subwords equal to $xyxy$
without reference. Besides that,~\Ref{a1} implies
$(xayx)y=-x^2ayy-a(yx^2y)=-x^2ay^2+ay^2x^2+ax^2y^2$, and
$xa(yxy)=-(xay^2x)-(xax)y^2=x^2ay^2+ay^2x^2+x^2ay^2+ax^2y^2$.
Hence
\begin{eq}\label{f2}
x^2ay^2=0, \mbox{ if } p\neq3.
\end{eq}
By separate linearization of \Ref{f2} with respect to $x$ and with
respect to $y$, we obtain
\begin{eq}\label{f3}x^2abc+x^2acb=0,abcx^2+bacx^2=0,\mbox{ if }
p\neq3.
\end{eq}
Applying~\Ref{a1}, we get $(xux)vx=-x^2uvx-ux^2vx$,
$xu(xvx)=-xux^2v-xuvx^2=x^2uxv+x^2uvx$. Hence
\begin{eq}\label{Azzz}
-2x^2uvx=x^2uxv+ux^2vx.
\end{eq}

\subsection{The case of $N_{3,2}$ in characteristic different from
$3$}
\begin{predl}\label{A2}
If $p\neq3$, then $C(3,2,K)=6$.
\end{predl}
\begin{proof} Applying \Ref{a} and \Ref{f2} to $x^2y^2xy$, we
obtain that $p\neq3$ implies $x^2y^2xy=0$. Statement~\ref{prop3}
concludes the proof.
\end{proof}

\begin{predl}\label{prop3}
$x_1^2x_2^2x_1\neq0$ for each $p$.
\end{predl}
\begin{proof} Let us find a solution for $\sys_{(3,2)}$ for which
$x_{x_1^2x_2^2x_1}\neq0$. Let $x_{x_1^2x_2^2x_1}=1$,
$x_{x_1x_2^2x_1^2}=-1$, $x_{x_1^2x_2x_1x_2}=-1$,
$x_{x_2x_1x_2x_1^2}=1$, $x_{x_1x_2x_1^2x_2}=1$,
$x_{x_2x_1^2x_2x_1}=-1$, $x_{x_1x_2x_1x_2x_1}=0$ and $x_u=0$ for
any other word $u$ of multidegree $(3,2)$. It is easy to see that
this is indeed a solution for every equation from $\sys_{(3,2)}$.
\end{proof}

\subsection{The case of characteristic equal to $0$ or greater than $3$}
\begin{prop}
If $p=0$ or $p>3$, $d>1$, then $C(3,d,K)=6$.
\end{prop}
\begin{proof} Equality $x_1x_2\cdots x_5=0$ implies $x^2y^2x=0$,
which is a contradiction to Statement~\thRef{A2}.

Applying \Ref{f3}, we get $x^2\cdot a\cdot b\cdot cd=-x^2\cdot
a\cdot cd\cdot b$ and $(x^2abc)d=-(x^2acb)d=-x^2\cdot ac\cdot
b\cdot d=x^2\cdot ac\cdot d\cdot b$. Hence
\begin{eq}\label{5.1}x^2abcd=0, abcdx^2=0.\end{eq}
It follows that $x^2y^2ab=0,abx^2y^2=0$. Linearization of these
identities with respect to $x$ gives
\begin{eq}\label{5.24}
abx^2cd+bax^2cd=0, abx^2cd+abx^2dc=0.
\end{eq}
Further,~\Ref{f3} and \Ref{5.24} imply $T_3(x^2ab,c,d) =
cx^2abd+dx^2abc = 0$ and $T_3(x^2,a,b)cd = ax^2bcd+bx^2acd=0$.
These two identities together with \Ref{f3} imply that
\begin{eq}\label{5.4}
a_1x^2a_2a_3a_4=\sign\si\cdot
a_{\si(1)}x^2a_{\si(2)}a_{\si(3)}a_{\si(4)}\mbox{, }\si\in S_4.
\end{eq}
Note that the same is true for $abcx^2d$. Further, $T_3(ax^2b,c,d)
= cax^2bd+dax^2bc=0$ (see~\Ref{5.4}). This identity together
with~\Ref{5.24} imply
\begin{eq}\label{5.5}
a_1a_2x^2a_3a_4=\sign\si\cdot
a_{\si(1)}a_{\si(2)}x^2a_{\si(3)}a_{\si(4)}\mbox{, }\si\in S_4.
\end{eq}
Let $A=ax^2bcd$, $B=abx^2cd$, $C=abcx^2d$. Then \Ref{5.1},
\Ref{5.4}, \Ref{5.5} imply $T_3(x^2,a,bc)d=A+B+2C=0$,
$dT_3(x^2,a,bc) = -2A-B-C=0$, $T_3(x^2a,bc,d)=-A+B-2C=0$.
Hence $A=B=C=0$. Linearization of~\Ref{5.1} and $A=B=C=0$ gives
$a_1\cdots a_6=\sign\si\cdot a_{\si(1)}\cdots a_{\si(6)}$, $\si\in
S_6$. It follows that $T_3(ab,cd,ef)=0$ implies $abcdef=0$.
\end{proof}

\subsection{The case of characteristic $2$}\label{case_p_2}
\begin{prop}
If $p=2$, then $C(3,d,K)=\left\{
\begin{array}{lcl}
d+3&,& d\geq3\\
6&,& d=2.\\
\end{array}
\right. $
\end{prop}
\begin{proof}
If $d=2$, see Statement~\thRef{A2}.

A word of multidegree $\Lambda=(\la_1,\ldots,\la_d)$, where
$\la_i>1$ $(i=\overline{1,3})$, is equal to $0$ by \Ref{f2}. Let
us show that a word $u$ of multid¥gree $\Lambda$, where
$\la_1=3,\la_2>1$, $\la_3>0$, is equal to $0$. Applying
Statement~\ref{canon}, we represent $u$ as a sum of words
containing subwords $x^2y^2xz$ (see \Ref{f2},\Ref{a}). But
$x^2\cdot y^2\cdot x\cdot z=/\mbox{see }
\Ref{f3}/=x^2y^2zx=/\mbox{see } \Ref{a}/=xy^2zx^2= /\mbox{see }
\Ref{f2}/=0$.

Let $d=3$. We have $x^2y^2x\neq0$ (by Statement~\ref{prop3}), and
all words of degree $6$ are equal to $0$. Hence, the nilpotency
degree is equal to $6$.

Let $d\geq4$. The longest words which can be non-zero are words of
multidegrees $(2,2,1,\ldots,1)$ and $\Theta=(3,1,1,\ldots,1)$.
Below we prove the existence of a non-zero word of multidegree
$\Theta$.
\end{proof}

\begin{predl}\label{predl_p2_1}
$x_1^2x_2\cdots x_dx_1\neq0$, where $d\geq2$.
\end{predl}
\begin{proof}
Let $V=x_1^2x_2\cdots x_{d}x_1$, $d\geq2$.  First let us show that
$V\neq0$ when $d\geq4$, which will imply that if $x_1^2x_2x_1=0$
or $x_1^2x_2x_3x_1=0$, then the substitution $x_2\rightarrow
x_2x_3x_4$ or $x_3\rightarrow x_3x_4$, respectively, leads to
required contradiction.

Let $d\geq4$. There exists a solution for ${\sys}_{\Theta}$ for
which $V\neq0$, namely, take $\{x_u=F(u)| \mdeg(u)=\Theta\}$,
where $F(u)$ is equal to the number of all subwords $x_1^2$ in
$u$. For example, if $\deg_{x_1}(u)=\deg_{x_1}(v)=0$, then
$F(ux_1^3v)=0$. Let
$$
F(v,w)=\left\{
\begin{array}{ccl}
1&,&\mbox{if } v=v'x_1, w=x_1w'\\
0&,&\mbox{otherwise }\\
\end{array}
\right..
$$
Here subwords $v'$, $u'$ can be empty. We have $F(v_1\cdots
v_l)=\sum\nolimits_{i=1}^{l} F(v_i) + \sum\nolimits_{i=1}^{l-1}
F(v_i,v_{i+1})$. Hence $g_1T_3(f_1,f_2,f_3)g_2|_{\{F(u)\}}=0$. Let
$g_1T_2(x_1,f)g_2=0$ and $g_1T_1(f)g_2=0$ be equations from
${\sys}_{\Theta}$. As one can see $g_1T_2(x_1,f)g_2|_{\{F(u)\}} =
F(g_1)+F(f)+F(g_2)+F(g_1,f)+F(f,g_2)$. We have
$\deg_{x_1}(g_1fg_2)=1$, so
$F(g_1)=F(f)=F(g_2)=F(g_1,f)=F(f,g_2)=0$. Hence,
$g_1T_2(x_1,f)g_2|_{\{F(u)\}}=0$. Clearly,
$g_1T_1(f)g_2|_{\{F(u)\}}=0$. So $\{x_u=F(u)\}$ is a solution for
${\sys}_{\Theta}$.
\end{proof}

\begin{note}\label{note_p2}
One can show that there exist non-zero words of multidegree
$(2,2,1,\ldots,1)$ (namely, $x_1^2x_2^2x_3\cdots x_d\neq0$, where
$d\geq2$).
\end{note}



\subsection{The case of characteristic $3$}\label{case_p_3}
\begin{prop}\label{prop6}
If $p=3$, then $C(3,d,K)=\left\{
\begin{array}{lcl}
3d+1&,& d\ {\rm is\ even}\\
3d\;{\rm or}\; 3d+1 &,& d\ {\rm is\ odd}.\\
\end{array} \right.
$
\end{prop}

The proof will follow from Statement~\ref{predl3} and
Corollary~\ref{predl3_1}.

\begin{predl}\label{predl2}
Let $\sum\al_iu_i=0$ be a homogeneous identity of degree $1$ or
$2$ in $x_k$, $k\in\overline{1,d}$, which contains some other
letters. Then the result of substitution $x_k=1$ in
$\sum\al_iu_i=0$ is an identity.
\end{predl}
\begin{proof}
Let $M$ be the set of all identities from system $\sys$ of degree
$1$ or $2$ in $x_k$. The identity $\sum\al_iu_i=0$ is a
consequence of identities from $M$. Set $M$ does not contain
identities $g_1T_1(f)g_2=0$, where $\deg_{x_k}(f)\neq0$. Hence the
result of substitution $x_k=1$ in any identity from $M$ is an
identity.
\end{proof}

Consider a word $x^2ux$, where $\deg_x(u)=0$. Replacing $x$ with
$x+y$, where $\deg_y(u)=0$, and taking the homogeneous component
of degree $1$ in $x$ and $2$ in $y$,  we get $y^2ux+yxuy+xyuy$.
Substitution $y=1$ gives $ux-xu$. This reasoning shows that linear
function $\Pi_{x}(v_1x^2uxv_2)=v_1uxv_2-v_1xuv_2$, where $v_1$,
$v_2$ are any words, is defined correctly on all homogeneous
components of $N_{3,d}$ of degree $3$ in $x$. Let
$W_{xy}=x^2y^2xy$. We will shorten $W_{x_ix_j}$ to $W_{ij}$, and
$\Pi_{x_i}$ to $\Pi_i$. We have
\begin{eq}\label{q1}
\Pi_{i}\Pi_{j}(W_{ij})=x_ix_j-x_jx_i.
\end{eq}
The element $u=x_{\pi(1)}\cdots x_{\pi(t)}\in K\LA
x_1,\ldots,x_t\RA$, $\pi\in S_t$, is called even if permutation
$\pi$ is even, and odd otherwise. Define $\sign u=1$ for even $u$
and $\sign u=-1$ for odd $u$. Denote by $|u|$ the length of $u$.
\begin{lemma}\label{lemma_q4} If $|v_1|$ and $|v_2|$ are both odd or both
even, then
$$\sign uv_1u'v_2u''=(-1)^{|v_1|\cdot|v_2|}\sign
uv_2u'v_1u'',$$%
where words $u$, $u'$, $u''$ can be empty.
\end{lemma}
\begin{proof}
The statement follows from $\sign
uw_1w_2u''=(-1)^{|w_1|\cdot|w_2|}\sign uw_2w_1u''$, where words
$u$, $u''$ can be empty.
\end{proof}
\begin{predl}\label{predl3} The word $w_{2k}=W_{12}W_{34}\cdots
W_{2k-1,2k}$ is not equal to $0$, if $k\geq1$.
\end{predl}
\begin{proof}
Assume that, on the contrary, $w_{2k}=0$. Then let
$$h_{2k}=\Pi_{1}\Pi_{2}\cdots\Pi_{2k-1}\Pi_{2k}(w_{2k})=(x_1x_2-x_2x_1)\cdots(x_
{2k-1}x_{2k}-x_{2k}x_{2k-1})=0.$$ If a word $u\in K\LA
x_1,\ldots,x_t\RA$ of multidegree $1^{(2t)}$ is even, let
$N_{+}(u)=1$ and $N_{-}(u)=0$; if it is odd, let $N_{+}(u)=0$ and
$N_{-}(u)=1$. Let us show that for every equation
$h=g_1T_3(f_1,f_2,f_3)g_2=0$ from ${\sys}_{\Lambda}$ (where
$\Lambda=1^{(2k)}$) it is true that
\begin{eq}\label{q2}
h|_{\{N_{+}(u)\}}=h|_{\{N_{-}(u)\}}=0.
\end{eq}
It is enough to consider equations with $g_1=g_2=1$, because for
words $uv_1$, $uv_2$ of multidegree $\Lambda$ if $\sign v_1=\sign
v_2$, then $\sign uv_1=\sign uv_2$. There are two possibilities:

{$1)$} Among $f_1$, $f_2$, $f_3$ there are two words of odd
length, for example, $f_2$ and $f_3$. Then by
Lemma~\ref{lemma_q4}, $\sign f_1f_2f_3=-\sign f_1f_3f_2$, $\sign
f_2f_1f_3=-\sign f_3f_1f_2$, $\sign f_2f_3f_1=-\sign f_3f_2f_1$.
Hence~\Ref{q2} is true.

{$2)$} Among $f_1$, $f_2$, $f_3$ there are two words of even
length, for example, $f_2$ and $f_3$. Then by
Lemma~\ref{lemma_q4}, words $f_1f_2f_3$, $f_1f_3f_2$, $f_2f_1f_3$,
$f_2f_3f_1$, $f_3f_1f_2$, $f_3f_2f_1$ are all even or all odd,
hence~\Ref{q2} is true.

Prove by induction on $t$ that
${h_{2t}}|_{\{N_{+}(u)\}}=(-1)^{t+1}$,
${h_{2t}}|_{\{N_{-}(u)\}}=(-1)^{t}$. For $h_2$ it is obvious.
Since $h_{2t}=h_{2(t-1)}(x_{2t-1}x_{2t}-x_{2t}x_{2t-1})$, we have
${h_{2t}}|_{\{N_{+}(u)\}} =h_{2(t-1)}|_{\{N_{+}(u)\}}
-{h_{2(t-1)}}|_{\{N_{-}(u)\}}$ and ${h_{2t}}|_{\{N_{-}(u)\}} =
{h_{2(t-1)}}|_{\{N_{-}(u)\}} -{h_{2(t-1)}}|_{\{N_{+}(u)\}}$. By
induction hypothesis, we get what is required.

We found a solution for $\sys_{\Lambda}$ on which $h_{2k}$ is not
equal to zero, which is a contradiction.
\end{proof}

\begin{concl}\label{predl3_1}
The word $W_{12}W_{34}\cdots W_{2k-1,2k}x_{2k+1}^2$ is not equal
to zero if $k>0$.
\end{concl}
\begin{proof}
The proof follows from Statements~\ref{predl2} and~\ref{predl3}.
\end{proof}


In Section 3 we will need Statement~\ref{predl2.p3}, which is
formulated below. We have $0=(xy)^3=(xyx)(yxy)= /\mbox{see }
\Ref{a1}/=(x^2y+yx^2)(y^2x+xy^2)=/\mbox{see
}\Ref{a}/=-x^2y^2xy-y^2x^2yx$. Hence
\begin{eq}\label{1p3}
W_{xy}=-W_{yx}.
\end{eq}
We will shorten $W_{xy}$ to $W$.
\begin{predl}\label{predl1}
Any word of degree $3$ with respect to $x$ and $y$ is equal to
$\sum \alpha_iu_iWw_i$, where subwords $u_i$, $w_i$ do not contain
$x$ and $y$.
\end{predl}
\begin{proof} By Statement~\ref{canon},
it is enough to consider canonical words.  For words of
multidegree $(3,3)$, Statement~\ref{predl1} follows
from~\Ref{1p3}. Let us prove it for words of multidegree
$(3,3,1)$.
$$\begin{array}{rcl}
T_3(x^2y^2,a,xy)=0&\Rightarrow &x^2y^2axy=-aW-Wa.\\
T_2(xy,ayx)=0&\Rightarrow &x^2y^2ayx=aW-Wa.\\
T_3(x^2,ay^2x,y)=0&\Rightarrow &x^2ay^2xy=aW.\\
T_3(y^2,yx^2a,x)=0&\Rightarrow &x^2y^2xay = Wa.\\
\end{array}$$

Consider identities of multidegree $(3,3,1,1)$. Identity
$T_3(x^2a,xy^2,by)=0$ implies $x^2axy^2by=abW+Wab-baW-Wba$. The
latter identity, together with $T_3(a,b,W)=0$, implies
$$x^2axy^2by=-abW-Wab+aWb+bWa.$$
We apply~\Ref{Azzz} to subwords which are put into parentheses.
$$\begin{array}{l}
x^2ay^2bxy=(x^2\cdot a\cdot y^2b\cdot x)y=abW-Wab+bWa.\\
x^2ay^2xby=(x^2\cdot a\cdot y^2\cdot x)by=-abW-Wab-aWb+bWa.\\
x^2y^2axby=x^2(y^2\cdot ax\cdot b\cdot y)=-abW+Wab+bWa.\\
x^2ay^2byx=(x^2\cdot a\cdot y^2by\cdot x)=-Wab+bWa.\\
x^2y^2aybx=(x^2\cdot y^2ay\cdot b\cdot x)=abW-bWa.\\
\end{array}$$

Likewise we obtain identities of multidegree $(3,3,1,1,1)$:
$$\begin{array}{l}
x^2ay^2bxcy=(x^2\cdot a\cdot y^2b\cdot x)cy=abcW+acWb+bcWa-aWbc-Wabc.\\
x^2ay^2bycx=(x^2\cdot ay^2\cdot byc\cdot
x)=abcW-abWc-acWb+aWbc+bWac-Wabc.\\
T_3(x^2axby^2,c,y)=0\Rightarrow\\
x^2axby^2cy=cabW-abWc+caWb-cbWa-aWcb+bWac-cWab+Wcab.
\end{array}$$

Modulo the identities which we obtained, each word is equivalent
to an element of the required form.
\end{proof}

\begin{concl}\label{predl1.1}
Let $w$ have degree $3$ in $x$ and $y$. Then the result of
substitution $\{x\rightarrow y,\;y\rightarrow x\}$ is $-w$.
\end{concl}
\begin{proof} Let us denote the result of substitution by $u$.
Statement~\ref{predl1} implies $w=\sum\alpha_iu_iW_{xy}w_i=/{\rm
see~\Ref{1p3}}/=-\sum\alpha_iu_iW_{yx}w_i=-u$.
\end{proof}

By Corollary~\ref{predl1.1} we have for any word $u$
\begin{eq}\label{0www}
W_{ij}uW_{kl}=W_{kl}uW_{ij}.
\end{eq}
Thus every permutation of subwords of the form $W_{ij}$, $i\neq
j$, does not change a word. So with abuse of notation we denote
all such words by one symbol $W$, i.e. $uWv$ equals $uW_{ij}v$ for
some $i\neq j$ such that letters $x_i,x_j$ are not contained in
$u,v$.
\begin{concl}\label{predl1.2} Every word of multidegree $3^{(2k)}$
equals $\alpha W^k$.
\end{concl}
\begin{proof}
See Statement~\ref{predl1}.
\end{proof}
\begin{concl}\label{predl1.3}
Every word of multidegree $3^{(2k+1)}$ equals $\alpha
x^2WxW^{k-1}$, $k>0$.
\end{concl}
\begin{proof}
Identity $T_2(W,x)=0$ implies $W^2x+WxW+xW^2=0$. Multiplying the
latter identity first from the left and then from the right by
$x^2$, we get
\begin{eq}\label{bwww}
x^2W^2x=-x^2WxW, \;\;x^2W^2x=-Wx^2Wx
\end{eq}
(see~\Ref{a}). Thus
\begin{eq}\label{bwww1}
Wx^2Wx=x^2WxW.
\end{eq}
Identities $T_2(W,xW)=xW^3+W^2xW+WxW^2=0$ and
$T_2(W,Wx)=W^3x+W^2xW+WxW^2=0$ imply that
\begin{eq}\label{cwww}
xW^3=W^3x.
\end{eq}
Let $r\in\{0,1,2\}$, $s\geq0$. Identities~\Ref{cwww}, \Ref{bwww}
imply $x^2W^{3s+r}x=x^2W^rxW^{3s}=r\cdot x^2WxW^{3s+r-1}$, since
$p=3$. Owing to~\Ref{bwww1}, we have
\begin{eq}\label{dwww}
W^ix^2W^lxW^j=l\cdot x^2WxW^{i+j+l-1},\;{\rm where}\;\;i,j,k\geq0.
\end{eq}
This formula and Statement~\ref{predl1} conclude the proof.
\end{proof}
\begin{predl}\label{predl2.p3}
If $\mdeg(uv)=3^{(d)}$, $d=2k$ or $d=6m+1$ ($k,m>0$), then
$uv=vu$.
\end{predl}
\begin{proof} Owing to Statement~\ref{canon}, we may assume that the words $u$,
$v$ are canonical. We will prove the Statement by 'decreasing'
induction on $s$, where $s$ is the number of subwords of the form
$W$ in words $u$ and $v$.

{\bf Induction base.} Let $d=2k$. If $s=k$, then $u=W^l$,
$v=W^{k-l}$ $(0<l<k)$, and Statement~\ref{predl2.p3} follows
from~\Ref{0www}.

Let $d=6m+1$. If $s=3k$, then $uv=W^ix^2W^lxW^j$, where
$i+j+l=3m$, and both subwords $u$ and $v$ are products of elements
of the set $\{x^2,x,W,\ldots,W\}$. Consider all possibilities:

1) $uv=W^ix^2W^{l_1}\cdot W^{l_2}xW^j$. Identities~\Ref{dwww}
and~\Ref{a} imply $uv=(l_1+l_2)x^2WxW^{3m-1}$,
$vu=-(j+i)x^2WxW^{3m-1}$. Since $i+j+l_1+l_2=3m$,  words $u$ and
$v$ commute.

2) $uv=W^ix^2W^lxW^{j_1}\cdot W^{j_2}$.

3) $uv=W^{i_1}\cdot W^{i_2}x^2W^lxW^j$.

The last two cases are similar to the first one.

{\bf Induction step.} Assume that $x$ and $y$ are not contained in
any of the subwords $W$ of words $u$ and $v$. Up to change of
notations, all possibilities can be reduced, by means of~\Ref{a},
to the following:

1) $u$ contains $x^2$, $x$, $y^2$, $y$;  $v$ does not contain
letters $x$, $y$.

2) $u=ax^2by^2c$, $v=dxeyf$.

3) $u=ax^2by^2c$, $v=dyexf$.

4) $u=ax^2bxc$, $v=dy^2eyf$.

5) $u=ax^2bxcy^2d$, $v=eyf$.

6) $u=ax^2by^2cxd$, $v=eyf$.

7) $u=ay^2bx^2cxd$, $v=eyf$.

Here words $a,\ldots,f$ can be empty. Consider these cases:

1) is obvious.

2) Identity $uv=vu$ is equivalent to $a(x^2by^2c\cdot
dxey)f=d(xeyf\cdot ax^2by^2)c$. Applying Statement~\ref{predl1} to
the subwords in parentheses, we can see that the previous identity
is equivalent to $abcdeWf+abeWcdf+acdeWbf-abWcdef-aWbcdef =
defabWc+debWfac+ dfabWec- deWfabc-dWefabc$. By induction
hypothesis, $-abcdeWf+ abeWcdf+acdeWbf+ abWcdef-aWbcdef -acdebWf
-abWecdf +abcdWef = 0$. Changing notations $b\rightarrow a$,
$cd\rightarrow b$, $e\rightarrow c$, we can see that the latter
identity follows from $-abcW-bcaW+abWc+acWb+bcWa+aWbc -aWcb-Wabc =
0$. The last identity was verified by means of a computer
programme, which was written by means of Borland C++ Builder
(version 5.0). The programme is available upon request from the
author.

The rest of the possibilities can be treated likewise.
%
%
%
%
\end{proof}

\section{Matrix invariant algebra}\label{chapter_matr}
\subsection{Auxiliary results}\label{technik}
Similarly to the definition of $R_d$ in terms of projective
limits, let
$\ov{R_d}=\bigoplus_{r\geq0}\projlim_n\ov{R_{n,d}}(r)$, or,
equivalently, $\ov{R_d}=R_d/(R_d^{+})^2$. The algebras $R_{n,d}$,
$\ov{R_{n,d}}$, $R_d$, $K\mbox{-}\alg\{X_1,\ldots,X_d\}\subset
C_{n,d}$ possess the natural $N_0$-grading by degrees and
$N_0^d$-grading by multidegrees, for which we use the same
notations as in Introduction. If elements $r_1,r_2$ of $R_{n,d}$
(of $R_d$, respectively) are equal modulo the ideal
$(R_{n,d}^{+})^2$ ($(R_d^{+})^2$, respectively), we write
$r_1\equiv r_2$. Since the ideal $(R_{n,d}^{+})^2$ is homogeneous
with respect to the $N_0^d$-grading, one can see that for every
equality of the form $r\equiv0$, $r\in R_{n,d}$, and
$N_0^d$-homogeneous component $r^\prime$ of $r$, $r^\prime\equiv0$
is also true. As in Section~2, monomials in the generic matrices
$X_i\in C_{n,d}$ are called words, and $X_i$ --- letters. The same
terminology is used for elements of $C_d$. By letters $U,V,W$,
possibly, with indices, we denote non-empty words in the generic
matrices, if we do not explicitly write otherwise.

\begin{lemma}\label{lemma2}
$$\frac{C_{n,d}}{\mbox{\rm id}\{ R_{n,d}^{+}\}} \simeq N_{n,d}.$$
\end{lemma}
\textbf{Proof}. As we mentioned in Introduction,
$$C_{n,d}\simeq \frac{C_d}{J_{n,d}}\simeq %
\frac{R_d\LA x_1,\ldots,x_d\RA}{\mbox{\rm id}\{\chi_k(f)| k\geq
n,\;f\in R_d\LA x_1,\ldots,x_d\RA^{\#}\}}.$$ %

Let $f\in C_d^{\#}$ and $f=f'+f''$, where $f'\in K\LA
x_1,\ldots,x_d\RA^{\#}$, $f''\in R_{d}^{+}C_d^{\#}$. Then
$\chi_k(f)\in C_d$ is equal to $f'^k$ modulo the ideal
$\id{R_{d}^{+}}\triangleleft C_d$, because $\sigma_j(g)\in
R_d^{+}$ for every $g\in C_{d}^{\#}$, $j>0$. Thus the ideal
$J_{n,d}$ is equal to $\mbox{\rm id}\{f^n| f\in K\LA
x_1,\ldots,x_d\RA^{\#}\}\triangleleft C_d$ modulo the ideal
$\mbox{\rm id}\{R_{d}^{+}\}$. It is easy to see that the preimage
of the ideal $\mbox{\rm id}\{R_{n,d}^{+}\}\triangleleft C_{n,d}$
in $C_d$ is equal to $\mbox{\rm id}\{R_{d}^{+}\}+J_{n,d}$. By the
Theorem on Homomorphism and the two preceding remarks, we have

$$
\frac{C_{n,d}}{\mbox{\rm id}\{ R_{n,d}^{+}\}} \simeq
\left.{\frac{C_d}{J_{n,d}}}\right/\frac{\mbox{\rm id}\{
R_d^{+}\}+J_{n,d}}{J_{n,d}}\simeq%
\frac{C_d}{\mbox{\rm id}\{ R_d^{+}\}+J_{n,d}}\simeq$$

$$
\simeq\frac{C_{d}}{\mbox{\rm id}\{ R_d^{+}\}+\mbox{\rm id}\{ f^n|
f\in K\LA x_1,\ldots,x_d\RA^{\#}\}}\simeq N_{n,d}. \bigtriangleup
$$

The image of $G\in C_{n,d}$ in $N_{n,d}$ we denote by $\ov{G}$. We
denote any triple of pairwise distinct generic matrices by
$X,Y,Z$, and their images we denote by  $x,y,z$. We assume that
$\si_k(f)$, where $f\in C_d^{\#}$, is an element of $R_{n,d}$,
unless it is stated otherwise.

We will use the fact that $\tr(XY)$ is a nongenerate bilinear
form, namely, if $\tr(GX_d)=0$, where $G\in M_n(K_{n,d-1})$, then
$G=0$.

\begin{lemma}\label{cr_1}
Suppose that $G\in K\mbox{-}\alg\{X_1,\ldots,X_d\}^{\#}$. Then

1) If $\ov{G}=0$ in $N_{n,d}$ then $\si_k(GX)$ is decomposable,
where $k>0$.

2) If $G$ does not contain $X$ and $\tr(GX)$ is decomposable, then
$\ov{G}=0$ in $N_{n,d}$.
\end{lemma}
\begin{proof}
1) Owing to Lemma~\ref{lemma2}, the identity $\ov{G}=0$ implies
$G=\sum_i r_iU_i$, where $r_i\in R_{n,d}^{+}$, $\deg(U_i)\geq0$.
Thus $\si_k(GX)=\si_k(\sum_{i} r_iU_iX)$ is decomposable
(see~\Ref{amiss}).

2) Let $\tr(GX)$ be decomposable. Since $R_{n,d}$ is homogeneous,
we have $\tr(GX)=\sum\tr(U_iX)r_i$, where words $U_i$ may be
empty, $\deg_X(U_i)=0$ and $r\in R_{n,d}^{+}$. Hence $\tr((G-\sum
r_iU_i)X)=0$. Since the trace form is nongenerate, we have $G=\sum
r_iU_i\in \mbox{id}\{R_{n,d}^{+}\}$. Therefore $\ov{G}=0.$
\end{proof}

Further, we assume $n=3$. Word $U$ is called  {\it canonical}, if
$\overline{U}\in N_{3,d}$ is canonical.

\begin{lemma}\label{can}
For every non-empty word $U$ there exist decompositions
$$\tr(U)\equiv\sum\alpha_i\tr(W_i),\;\;\Si(U)\equiv\sum\beta_j\Si(U_j)+\sum
\gamma_l\tr(V_l),$$ where $W_i,U_j,V_l$ are canonical words.
Moreover, homogeneity of ideal $(R_{3,d}^{+})^2$ implies that
multidegrees of $\tr(U)$ and $\tr(W_i)$ are all equal, and also
that multidegrees of $\si_2(U)$, $\si_2(U_j)$ and $\tr(V_l)$ are
all equal.
\end{lemma}
\begin{proof} First we will prove the Lemma for the trace. Owing to Lemma~\ref{cr_1}
 and linearity of the trace, one can prove formulas analogous to~\Ref{a1} and~\Ref{a},
 namely, for $i=\ov{1,d}$ we have $\tr(V_1X_iVX_iV_2) \equiv -\tr(V_1X_i^2VV_2) -
 \tr(V_1VX_i^2V_2)$,
\begin{eq}\label{can_id2}
\tr(V_1X_iVX_i^2V_2)\equiv-\tr(V_1X_i^2VX_iV_2).
\end{eq}
Here one of words $V_1$, $V_2$ can be empty. Hence the Lemma is
proved for the trace.

The proof for $\Si(U)$ is similar, except that instead of the
trace linearity we apply consequence of Amitsur's formula for
$\Si$: $\Si(V_1+V_2)\equiv\Si(V_1)+\Si(V_2)-\tr(V_1V_2)$, and then
we apply the proved part of the Lemma to $\tr(V_1V_2)$.
\end{proof}

\begin{lemma}\label{cr_2}
Suppose that $G\in K\mbox{-}\alg\{X_1,\ldots,X_d\}$, $G$ does not
contain $X$. Then

1)If $\tr(GX^2)$ is decomposable then $gx+xg=0$ in $N_{3,d}$,
where $g=\ov{G}$.

2)In the case $p\neq2$, the converse is also valid.
\end{lemma}
\begin{proof}
1) Substituting $X+Y$ for $X$, where $X$, $Y$ are not contained in
$G$, and taking the homogeneous component of degree $1$ in $X$ and
in $Y$, we get $\tr(GXY)+\tr(GYX)\equiv0$. Hence
$\tr((GX+XG)Y)\equiv0$. Lemma~\ref{cr_1} concludes the proof.

2) By Lemma~\ref{cr_1}, we have $\tr((GX+XG)X)\equiv0$. Hence,
$2\tr(GX^2)\equiv0$.
\end{proof}

\begin{lemma}\label{cr_3}
Let $\deg_X(U)=\deg_X(V)=0$ and $u=\ov{U}$, $v=\ov{V}$.

1)If $\tr(X^2UXV)\equiv0$ then $ux^2v-2vx^2u-x^2uv-uvx^2=0$ in
$N_{3,d}$.

2)In the case $p\neq3$, the converse is also valid.
\end{lemma}
\begin{proof}
1) Substituting $X+Y$ for $X$ in $\tr(X^2UXV)\equiv0$, and taking
the homo\-geneous component of degree $2$ in $X$ and of degree $1$
in $Y$, we get $\tr(VX^2UY)+\tr(UXVXY)+\tr(XUXVY)\equiv0$. Here
words $U$, $V$ do not contain $X$, $Y$. By Lemma~\ref{cr_1}, we
have $vx^2u+uxvx+xuxv=0\;\;\mbox{in } N_{3,d}$. This, together
with identity ~\Ref{a1}, yields the required equality.

2) By Lemma~\ref{cr_1}, we have $\tr(UX^2VX)-2\tr(VX^2UX)\equiv0$.
The identity~\Ref{can_id2} gives $3\tr(X^2UXV)\equiv0$.
\end{proof}

Applying Amitsur's formula to $\si_2(U+V)$ and letting $V=U$,
obtain
\begin{eq}\label{q0}
2\Si(U)=\tr^2(U)-\tr(U^2).
\end{eq}
\begin{lemma}\label{predl_si}
$a$) $\Si(X)$, $\det(X)$ are indecomposable. In particular,
$D_{\det}(3,d,K)=3$ and $D(3,1,K)=3$.

$b$) $\tr(X^2)$ is decomposable $\Leftrightarrow$ $p=2$.

$c$) $\tr(X^3)$ is decomposable $\Leftrightarrow$ $p=3$.
\end{lemma}
\begin{proof}
$a)$ Let $\si_k(X)$ be decomposable ($k=2,3$). Then, by $D)$,
$\si_k(X)$ can be expressed in terms of $\si_l(X)$, $l<k$.
Substitution $X=\diag(x_1,x_2,x_3)$ yields a contradiction to the
fact that an elementary symmetric polynomial can not be expressed
in terms of other elementary symmetric polynomials.

Lemma~\ref{cr_1} implies $\tr(X^k)\equiv0$, $k>3$, and $D)$
implies $\si_2(X^k)\equiv0$, $k>1$.

$b)$ For $p=2$, by~\Ref{q0}, $\tr(X^2)\equiv0$. For $p\neq2$, if
$\tr(X^2)\equiv0$, then $2x=0$ in $N_{3,d}$ (see
Lemma~\ref{cr_2}), which is false.

$c)$ If $\tr(\chi_3(X))=0$, then $\tr(X^3)\equiv3\det(X)$.
\end{proof}


\subsection{The case of characteristic equal to $0$ or greater than $3$}
\begin{predl}
If $p=0$ or $p>3$, $d>1$, then $D(3,d,K)=6$.
\end{predl}
\begin{proof}
Element $\tr(X_1\cdots X_7)$ is decomposable by Lemma~\ref{cr_1}
and identity $x_1\cdots x_6=$~$0$ in $N_{3,d}$. Hence
$\tr(U)\equiv0$ for every $U$ with $\deg(U)>6$.

Let us prove that $\tr(X^2Y^2XY)$ is indecomposable. Assume that
it is decomposable. Letting $U=X^2$, $V=X$ and applying
Lemma~\ref{cr_3}, we obtain $x^2y^2x-2xy^2x^2=0$ in $N_{3,d}$.
Then $x^2y^2x=0$ (see~\Ref{a}). But this yields a contradiction
(see Statement~\ref{prop3}).

Formula~\Ref{q0} concludes the proof.
\end{proof}


\subsection{The case of characteristic equal to $2$}

\begin{prop}

If $p=2$, then $D(3,d,K)=\left\{
\begin{array}{ccl}
d+2&,&d\geq4\\
6&,&d=2\;{\rm or}\; 3\\
\end{array} \right..
$
\end{prop}

The proof is a consequence of the following two statements.

\begin{predl}\label{predl_p2_tr}
If $p=2$, then $D_{\tr}(3,d,K)=\left\{
\begin{array}{ccl}
d+2&,&d\geq4\\
6&,&d=2\;{\rm or}\; 3\\
\end{array} \right..
$
\end{predl}
\begin{proof} By Lemma~\ref{can}, it is sufficient to consider
$\tr(U)$, where $U$ is canonical.

First we point out some restrictions on the multidegree of an
indecomposable invariant, namely, $\tr(U)$ is decomposable if
multidegree of $U$ is equal to

$a)$ $\Delta_1=(2,2,2,i_4,\ldots,i_d)$, where $i_4+\cdots+i_d\geq
1$.

$b)$ $\Delta_2=(3,2,i_3,i_4,\ldots,i_d)$, where
$i_3+\cdots+i_d\geq 2$.

Let us prove it. Every canonical word of multidegree $\Delta_1$ is
equal to
\\ $U_1X_{\pi(1)}^2U_2X_{\pi(2)}^2U_3X_{\pi(3)}^2U_4$, $\pi\in
S_3$, where some words $U_1,\ldots,U_4$ (but not all of them) can
be empty. Formula~\Ref{f2} and Lemma~\ref{cr_1} yield
$\tr(X^2V_1Z^2V_2)\equiv0$. Hence, decomposability is established
for $a)$.

Let the multidegree of $U$ be equal to $(3,2,1,1,i_5,\ldots,i_d)$.
Then $\tr(U)$ is decomposable, because each word from $N_{3,d}$ of
multidegree $(3,2,1,j_4,\ldots,j_d)$ is equal to $0$ (see
Section~\ref{case_p_2}); then we apply Lemma~\ref{cr_1}, which
gives decomposability for $b)$.

Let $d=2$. Then $\tr(X^2Y^2XY)$ is a maximal indecomposable
element (otherwise $y^2x^2y=0$ in $N_{3,d}$, by Lemma~\ref{cr_3},
but this is a contradiction, by Statement~\ref{prop3}). Invariants
of greater degree are, evidently, decomposable.

Let $d=3$. Then $\tr(X^2Y^2Z^2)$ is a maximal indecomposable
element (otherwise Lemma~\ref{cr_2} implies $x^2y^2z+zx^2y^2=0$ in
$N_{3,d}$, thus $x^2y^2x=xx^2y^2=0$ in $N_{3,d}$, by~\Ref{a}, but
$x^2y^2x\neq0$ --- see Statement~\ref{prop3}). All invariants of
greater degree are decomposable by $b)$. Also note that
$\tr(X^2Y^2XZ)$ is indecomposable, because, assuming that it is
decomposable and letting $Z=Y$, we get that $\tr(X^2Y^2XY)$ is
decomposable.

Let $d\geq4$. Invariant $\tr(X_1^2X_2X_1X_3\cdots X_d)$ is a
maximal indecomposable element, because $x_1^2x_2x_1x_3\cdots
x_{d-1}\neq0$  in $N_{3,d}$ (see Statement~\ref{predl_p2_1} and
Lemma~\ref{cr_1}). All words of greater degree are decomposable by
$a)$ and $b)$. Note that $\tr(X_1^2X_2X_1X_3\cdots X_d)$ is
indecomposable (see Remark~\ref{note_p2} and Lemma~\ref{cr_1}).
\end{proof}

\begin{predl}\label{predl_p2_sigma}
If $p=2$, then $D_{\si_2}(3,d,K)=\left\{
\begin{array}{ccc}
6&,&d\geq3\\
4&,&d=2\\
\end{array} \right..
$
\end{predl}
\begin{proof}
Applying Amitsur's formula to $\si_4(u+v)=0$, where $u,v\in K\LA
x_1,\ldots,x_d\RA^{\#}$ are words, and considering the result
modulo the ideal $(R_{3,d}^{+})^2$, we obtain
$\si_2(UV)+\tr(U^3V)+\tr(V^3U)+\tr(U^2V^2)\equiv0$, where $U,V$
are non-empty words in the generic matrices. Since $\tr(U^3V)$ is
decomposable (see Lemma~\ref{cr_1}), we have
\begin{eq}\label{m_q2} \si_2(UV)\equiv \tr(U^2V^2).
\end{eq}
Letting $U=X^2$, we obtain $\Si(X^2V)\equiv0$ (see
Lemma~\ref{cr_1}).

Let $U=X_1$, $V=X_2\cdots X_d$. Then~\Ref{m_q2}, together with the
identity of $N_{3,d}$ $(x_2\cdots x_d)^2-x_d^2\cdots x_2^2=0$
(which is a consequence of~\Ref{1p1}), to which we apply Lemma
~\ref{cr_1}, yields $\Si(X_1\cdots X_d)\equiv \tr(X_1^2(X_2\cdots
X_d)^2)\equiv \tr(X_1^2X_d^2\cdots X_2^2)$. Element $\tr(X^2Y^2)$
is indecomposable (otherwise Lemma~\ref{cr_2} implies
$x^2y+yx^2=0$ in $N_{3,d}$, thus $x^2y^2x=-y^2x^2x=0$, but the
last identity contradicts Statement~\ref{prop3}). This reasoning
and Statement~\ref{predl_p2_tr} imply that $\Si(X_1\cdots X_d)$ is
decomposable iff $d\geq 4$. By Lemma~\ref{can}, the
statement~\ref{predl_p2_sigma} is proved.
\end{proof}

\subsection{The case of characteristic equal to $3$}
\begin{prop}\label{p3}
If $p=3$, then $D(3,d,K)= \left\{
\begin{array}{ccl}
3d&,&d\;{\rm \ is\ even}\\
3d-1&,&d\equiv 3,5\;(\mod 6) \\
3d-1 \;{\rm or}\; 3d&,& d\equiv 1\;(\mod 6). \\
\end{array} \right.$
\end{prop}

To prove the proposition, we need more detailed study of
identities of $\ov{R_{3,d}}$.

\begin{lemma}\label{pr3}
If $\sum\alpha_iu_i=0$ in $N_{3,d}$, where $u_i$ are words, then
$\sum\alpha_i\tr(U_i)\equiv0$ in $\ov{R_{3,d}}$, where
$\ov{U_i}=u_i$. Note that this identity is a consequence of
$$\tr(dT_1(a)e)\equiv0,\;\tr(dT_2(a,b)e)\equiv0,\;\tr(dT_3(a,b,c)e)\equiv0,$$
where words $a,b,c\in K\LA x_1,\ldots,x_d\RA^{\#}$, and words
$d,e\in K\LA x_1,\ldots,x_d\RA$.
\end{lemma}
\begin{proof}
Denote by $A,B,C,D,E$ words in the generic matrices, where $D$ and
$E$ can be empty. Owing to Lemmas~\ref{predl_si} and~\ref{cr_1},
$\tr(DA^3E)$ is decomposable. Linearization yields
$\tr(DT_2(A,B)E)\equiv0,\;\;\tr(DT_3(A,B,C)E)\equiv0$.

Hence, if $g=0$ is an identity from ${\sys}$, then $\tr(G)=0$
holds in $\ov{R_{3,d}}$, where $\ov{G}=g$ (see
Section~\ref{intro}). Since all identities of $N_{3,d}$ are
consequences of system ${\sys}$,  every identity of $N_{3,d}$ has
a counterpart in $R_{3,d}$.
\end{proof}
\begin{lemma}\label{m_lemma1}
All identities of the algebra
$\ov{R_{3,d}}=R_{3,d}/(R_{3,d}^{+})^2$ are consequences of

(a) $\tr(dT_1(a)e)\equiv0,\;\tr(dT_2(a,b)e)\equiv0$,
$\tr(dT_3(a,b,c)e)\equiv0$.

(b) $\tr(ab)\equiv\tr(ba)$.

(c) $\si_2(a)\equiv \tr(a^2)$.

(d) $\det(ab)\equiv0$.

(e) $\si_k(a)\equiv0$, $k>3$.

Here words $a,b,c\in K\LA x_1,\ldots,x_d\RA^{\#}$, and words
$d,e\in K\LA x_1,\ldots,x_d\RA$.
\end{lemma}
\begin{proof} For $I\in\{A,B,C,D,E\}$ let $(\ov{I})$ be the
identity obtained by factorization of $(I)$ modulo the ideal
$(R_{3,d}^{+})^2$. Denote by $(I_w)$, $(\ov{I_w})$, respectively,
those identities of type $(I)$, $(\ov{I})$, respectively, in which
$g,h\in K\LA x_1,\ldots,x_d\RA^{\#}$ are words. Throughout this
proof we denote by letters $u,v$, possibly with indices, non-empty
words from $K\LA x_1,\ldots,x_d\RA$. The ideal of relations of
$R_d$ is generated by $(A_w)$ and $(D_w)$ (see~\cite{Zub2}). Thus
for the proof it is sufficient to show that $(\ov{A_w})$,
$(\ov{D_w})$, $(\ov{E})$ can be deduced from $(a)$--$(e)$. For
$(\ov{A_w})$ it is obvious. Consider $(\ov{D_w})$:
$\si_k(u^t)\equiv \alpha_{k,t}\si_{kt}(u)$, where $k\geq1$,
$t\geq2$.

Let $k=1$, $t=2$. Since
 $\si_2(u)$ is indecomposable (Lemma~\ref{predl_si}), we have
 $\al_{1,2}=1$, and $(\ov{D_w})$ follows from $(c)$.

Let $k=1$, $t=3$. Since
 $\tr(u^3)$ is decomposable and $\det(u)$ is indecomposable
 (Lemma~\ref{predl_si}), we have $\al_{1,3}=0$, and $(\ov{D_w})$
 follows from $(a)$.

If $k=1$, $t\geq4$, then $(\ov{D_w})$ follows from $(a)$, $(e)$.

If $k=2$, then $(\ov{D_w})$ follows from $(c)$ and $(a)$, $(e)$.

If  $k=3$, then $(\ov{D_w})$ follows from $(d)$, $(e)$.

If  $k\geq4$, then $(\ov{D_w})$ follows from $(e)$.

For the proof of deducibility of  $(\ov{E})$ we need some
properties of identities of $\ov{R_{3,d}}$.

\begin{enumerate}
\item[1) ] If an identity of $\ov{R_{3,d}}$ $t(x_1,\ldots,x_d)\equiv0$
can be deduced from $(a)$--$(e)$, then the identity
$t(r_1u_1,\ldots,r_du_d)\equiv0$, where $r_i\in R_{3,d}$, can be
deduced from $(a)$--$(e)$.
\end{enumerate}

Let us prove $1)$. By homogeneity of $(a)$--$(e)$ we may assume
$t(x_1,\ldots,x_d)\equiv0$ to be $N_0^d$-homogeneous. Then the
identity $t(r_1u_1,\ldots,r_du_d)\equiv0$ has the form $r\cdot
t(u_1,\ldots,u_d)\equiv0$, $r\in R_{3,d}$. Clearly the latter
identity is a consequence of $(a)$--$(e)$.

\begin{enumerate}
\item [2)] Let $u_i$ be words such that $\deg_{x_1}(u_i)\in\{1,2\}$,
and let
\begin{eq}\label{to1}
\sum\al_i\tr(u_i)\equiv0
\end{eq}
be an identity of $\ov{R_{3,d}}$. Then~\Ref{to1} follows  from
$(a)$, $(b)$. In particular, if~\Ref{to1} is an identity of
$\ov{R_{3,d}}$ and $\deg(u_i)\neq3s$, then~\Ref{to1} follows from
$(a)$, $(b)$.
\end{enumerate}

Let us prove $2)$. Identity~\Ref{to1} can be assumed to be
homogeneous. Let $\deg_{x_1}(u_i)= 1$. Rewrite identity~\Ref{to1}
in the form $\sum\al_i\tr(v_ix_1)\equiv0$, where words $v_i$ can
be assumed to be non-empty. By Lemma~\ref{cr_1}, $\sum \al_iv_i=0$
in $N_{3,d}$. But then $\sum \al_i v_ix_1=0$ in $N_{3,d}$, and
Lemma~\ref{pr3} concludes the proof.

Let $\deg_{x_1}(u_i)=2$. Identity~\Ref{to1} can be deduced from an
identity $\sum\beta_i\tr(v_ix_1^2)\equiv0$ by $(a)$ (see
Lemma~\ref{can}), where words $v_i$ can be assumed to be
non-empty. By Lemma~\ref{cr_2}, we have $\sum
\beta_i(x_1v_i+v_ix_1)=0$ in $N_{3,d}$. Substituting $x_1^2$ for
$x_1$, applying Lemma~\ref{pr3} and using $(b)$, we get the
required.

\begin{enumerate}
\item[3)] Denote by $(a_{*})$, $(b_{*})$, $(c_{*})$, $(d_{*})$
identities of $\ov{R_{3,d}}$ of the type $(a)$--$(d)$,
respectively, in which $a,b,c\in K\LA x_1,\ldots,x_d\RA^{\#}$ and
$d,e\in K\LA x_1,\ldots,x_d\RA$ (here $a,b,c,d,e$ are not
necessarily words). Thus $(a_{*})$--$(d_{*})$ follow from
$(a)$--$(d)$.
\end{enumerate}

Let us prove $3)$. Deducibility of $(a_{*})$, $(b_{*})$ from
$(a)$, $(b)$ is obvious. Owing to property $1)$, we can assume
that $a=\sum\nolimits_{i=1}^s x_i$, $b=\sum\nolimits_{i=s+1}^t
x_i$ in $(c_{*})$, $(d_{*})$.

Consider $(c_{*})$: $\si_2(\sum x_i)\equiv \tr((\sum x_i)^2)$.
Owing to $(c)$, identity $(c_{*})$ follows from some identity of
$\ov{R_{3,d}}$ of the type $\sum\beta_i\tr(v_i)\equiv0$, where
$v_i$ are words of degree $2$. The latter identity follows from
$(a)$, $(b)$ by property $2)$.

Consider $(d_{*})$: $\det((\sum_{i=1}^{s} x_i)(\sum_{i=s+1}^{t}
x_i))\equiv0$. By $(d)$, identity $(d_{*})$ follows from some
identity of $\ov{R_{3,d}}$ of the type $\sum
\gamma_l\tr(w_l)\equiv0$, where $w_l$ are products of words
$x_ix_j$ ($i=\ov{1,s}$, $j=\ov{s+1,t}$) and $\deg(w_l)=6$. Taking
homogeneous components, we may assume that this identity is
homogeneous of multidegree $\Theta$.

Let $\Theta\neq(3,3)$. Then all words $w_l$ have degree $1$ or $2$
in some letter $x_r$. Applying property $2)$, we conclude the
proof for this case.

Let $\Theta=(3,3)$. Then for each $l$ for some $r,q$ we have
$w_l=(x_rx_q)^3$, and the identity follows from $(a)$. Thus $3)$
is proved.

Now we can show that $(\ov{E})$:  $\si_k(h)\equiv0$, where
$k\geq4$, $h=\sum_{i=1}^{m} r_iu_i$, $r_i\in R_d$, follows from
$(a)$--$(e)$.

By property $1)$, we can assume $r_i=1$, $u_i=x_i$. Our proof is
by induction on $k\geq4$, and for a fixed $k$ --- by induction on
$m$.

{\bf Induction base}. Let us show that $(a)$--$(e)$ imply
$\si_k(x_1+x_2)\equiv0$. Owing to $(c)$--$(e)$, this identity is a
consequence of some identity of $\ov{R_{3,d}}$ of the form
$\sum\al_i\tr(v_i)\equiv0$.

If $k=4,5$, then $\deg(v_i)=4$ or $5$. The proof is concluded, by
property $2)$.

If $k=6$, then, by $(d)$ and $(e)$, the considered identity
follows from $-\si_2(x_1^2x_2) -\si_2(x_1x_2^2)+
\tr(x_1^2x_2^2x_1x_2)+ \tr(x_2^2x_1^2x_2x_1) \equiv0$. Identities
$(c)$ and $(a)$ imply $\si_2(x_1^2x_2)\equiv0$,
$\si_2(x_1x_2^2)\equiv0$. Identity
$\tr(x_1^2x_2^2x_1x_2)+\tr(x_2^2x_1^2x_2x_1)\equiv0$ follows from
$(a)$, by Lemma~\ref{pr3} applied to
$x_1^2x_2^2x_1x_2+x_2^2x_1^2x_2x_1=0$ in $N_{3,d}$
(see~\Ref{1p3}).

If $k\geq7$, then for every $i$ we have $\deg_{x_1}(v_i)>3$ or
$\deg_{x_2}(v_i)>3$, thus $v_i=0$ in $N_{3,d}$. Hence
$\tr(v_i)\equiv0$ follows from $(a)$ (see Lemma~\ref{pr3}).

{\bf Induction step}.
Consider identity of $\ov{R_{3,d}}$ $\si_k(x_1+x_2) = \si_k(x_1)
+\si_k(x_2) + \sum_j\al_{j}\si_{k_j}(u_j)\equiv0$, where $k_j<k$.
Let $g=\sum_{i=2}^m x_i$. The induction hypothesis yields
$\si_k(x_2|_{x_2\rightarrow g})\equiv0$,
$\si_{k_j}(u_j|_{x_2\rightarrow g})\equiv0$ $(k_j>3)$ follow from
$(a)$--$(e)$. Because $\si_k(x_1+x_2)\equiv0$ is a consequence of
$(a)$--$(e)$, we have $\sum_{k_j\leq3} \al_j\si_{k_j} (u_j)
\equiv0$ follows from $(a)$--$(d)$. Hence $\sum_{k_j\leq3}
\al_j\si_{k_j} (u_j|_{x_2\rightarrow g}) \equiv0$ follows from
$(a)$--$(d)$, by property $3)$. The lemma is proved.
\end{proof}

\begin{predl}\label{pr4}
Let $U_i$ be words of equal multidegree $\Theta$. Then
\begin{eq}\label{x1}
\sum\limits \alpha_i\tr(U_i)\equiv0
\end{eq}
is an identity of $\ov{R_{3,d}}$ if and only if system
${\sys_{\Theta}}$ and identities $uv=vu$, where $u,v\in K\LA
x_1,\ldots,x_d\RA^{\#}$ are words and $\mdeg(uv)=\Theta$, imply
that $\sum\limits \alpha_i\ov{U_i}=0$.
\end{predl}
\begin{proof}
$\Leftarrow$ Apply Lemma~\ref{pr3}.

$\Rightarrow$ By Lemma~\ref{m_lemma1}, identity~\Ref{x1} follows
from $(a)$--$(e)$. Identities~\Ref{x1}, $(a)$, $(b)$ do not
contain $\si_k(a)$, $k\geq2$, while identities~$(c)$--$(e)$
contain them. Every 'symbolic' element $\si_k(a)$, where $k\geq2$,
occurs in exactly one identity from $(c)$--$(e)$. Then the
derivation of~\Ref{x1} from~$(a)$--$(e)$ can be transformed into a
derivation of~\Ref{x1} from~$(a)$, $(b)$. The identities~\Ref{x1},
$(a)$, $(b)$ are homogeneous, thus for derivation of~\Ref{x1} we
only
need identities~$(a)$,~$(b)$ of multidegree $\Theta$. 
\end{proof}

Now we can prove Proposition~\ref{p3}.

\begin{proof}
By Lemma~\ref{can} and formula~\Ref{q0}, it is sufficient to
consider invariants of the form $\tr(U)$, where word  $U$ is
canonical. All words of the form $X_i^2X_j^2X_iX_j$, $i\neq j$,
are denoted by the same symbol $W$, and let $w=\ov{W}$ (see
Section~\ref{case_p_3} for details).

Let $d=2k$, $k>0$. In $N_{3,d}$ $w^k\neq0$, and also $uv=vu$,
where $\mdeg(uv)=3^{(2k)}$ (Statement~\ref{predl2.p3}).
Statement~\ref{pr4} yields the indecomposability of $\tr(W^{k})$.

Let $d=2k+1$, $k>0$. Invariant $\tr(X^2W^{k})$ is indecomposable,
because otherwise Lemma~\ref{cr_2} implies that $x_1w^k+w^kx_1=0$
in $N_{3,d}$. Substitution $x_1=1$ (see Statement~\ref{predl2})
yields $w^{k}=0$, which is a contradiction to
Statement~\ref{predl3}.

Let $d=6m+r$, $r\in\{3,5\}$, $m>0$. Let us show that if
$\mdeg(U)=3^{(d)}$, then $\tr(U)\equiv0$. For $u=\ov{U}$ we have
$u=\al v_d$, where $v_d=x^2wxw^{k-1}$, $d=2k+1$
(Corollary~\ref{predl1.3}). Hence $\tr(U)\equiv\al\tr(V_d)$, where
$\ov{V_d}=v_d$ (see Lemma~\ref{pr3}). If $r=3$, then~\Ref{cwww}
implies $\tr(V_d) \equiv \tr(X^2WW^{3m}X) \equiv0$. If $r=5$, then
it is easy to see that $\tr(V_d) =\tr(X^2WXW^{3m+1})=
\tr(XWW^{3m}X^2W)\equiv\tr(XWX^2W^{3m}W) \equiv /{\rm
see}~\Ref{a}/\equiv- \tr(X^2WXW^{3m+1})$. Hence $\tr(V_d)\equiv0$.
\end{proof}

\begin{center} { ACKNOWLEDGEMENTS} \end{center}

The author is grateful to A.N.Zubkov for helpful advices and
constant attention, to G.A.Bazhenova for help with the
translation. The author is also grateful to the referee whose
comments considerably improved the paper.

\end{document}